\documentstyle{amsppt}

\magnification=1200 \NoBlackBoxes \hsize=11.5cm \vsize=18.5cm

\def\inv{^{-1}}

\def\E{\Cal E}

\def\P{\Bbb P}

\def\O{\Cal O}

\def\Sym{\text{Sym}}

\def\1/2{\frac{1}{2}}

\def\I{\Cal I}

\def\Sec{\text{Sec}}

\def\im{\text{im}}

\def\2{{[2]}}
%\nologo
\topmatter
\title On the irreducibility of secant cones, and an
application to linear normality\endtitle
\author Angelo Lopez and Ziv Ran\endauthor
\address Universit\`a di Roma Tre\endaddress
\email lopez\@matrm3.mat.uniroma3.it\endemail
\address University of California, Riverside\endaddress

%\author xx\endauthor

\email ziv\@math.ucr.edu \endemail
%\author yy\endauthor
\rightheadtext {secant cones and linear normality}
\leftheadtext{Lopez and Ran} \abstract Given a smooth subvariety
of dimension $>\frac{2}{3}(r-1)$ in $\P^r$, we show that the
double locus (upstairs) of its generic projection to $\P^{r-1}$ is
irreducible. This implies a version of Zak's Linear Normality
theorem.\endabstract
\endtopmatter
\document
A classical, and recently revisited (cf. \cite {GP, L, Pi} and
references therein), method for studying the geometry of a
subvariety $Y$ in $\P^r$ is to project $Y$ generically to a
lower-dimensional projective space, for example so that $Y$ maps
birationally to a (singular) hypersurface
$\bar{Y}\subset\P^{m+1}.$ To make use of this method, it is
usually important to have precise control over the singularities
of $\bar{Y}$ and in particular over the entire singular (=double)
locus $D_Y$ of $\bar{Y}$ and its inverse image $C_Y$ in $Y$. As
the dimension of these is easily determined, a natural question
is: are $C_Y$ and $D_Y$ irreducible? This question plays an
important role, for instance, in Pinkham's work on regularity
bounds for surfaces \cite{Pi}.  The purpose of this note is to
show that this irreducibility holds provided the codimension of
$Y$ is sufficiently small compared to its dimension (see Theorems
1,2 and Corollary 3 below). As an application we give a proof of
Zak's linear normality theorem (in a slightly restricted range,
see Corollary 4 below). Indeed the results seem closely related as
our argument ultimately depends on having a bound on the dimension
of singular loci of hyperplane sections, manifested in the form of
the integer $\sigma(Y)$ (see Thm. 1 below), and it is Zak's
theorem on tangencies- also a principal ingredient in other proofs
of linear normality- that gives us good control over $\sigma(Y)$.

We begin with some definitions. Let $Y$ denote an irreducible $m-$
dimensional subvariety of $\P^r$. As usual, we mean by a
{\it{secant line}} of $Y$ a limit of lines in $\P^r$ spanned by
pairs of distinct points of $Y$. The union of all secant lines is
denoted by $\Sec(Y)$. $Y$ is said to be {\it{projectable}} if
$\Sec(Y)\subsetneq \P^r$. For any linear subspace $Q\subset \P^r$,
we let
$$\pi_Q:\P^r-Q\to\P^{r-\dim Q-1}$$
denote  the associated projection.\par For a nondegenerate
projective variety $Y$, let $\sigma(Y)$ denote the maximum
dimension of a subvariety $Z\subset Y_{\text{smooth}}$ such
that\par (i) $Z$ contains a generic point of some divisor on
$Y$;\par (ii) the tangent planes $T_yY$ for all $y\in Z$ are
contained in a fixed hyperplane $H$ (i.e. $Z$ is contained in the
singular locus of $Y\cap H$).\par Note that if $Y$ is nonsingular
in codimension 1 then assumption (i) above for a subvariety
$Z\subset Y$ already implies that $Z\cap Y_{\text{smooth}}$ is
dense in $Z$ . Zak's Tangency theorem (cf. \cite{F, Z1}) implies
that if $Y$ is smooth then
$$\sigma(Y)\leq r-m-1.\tag 1$$
\proclaim{Theorem 1} Let
$$Y\subset\P^r$$
be $m$-dimensional, normal, irreducible and non-projectable, and
let $Q\subset \P^r$ be a generic (resp. arbitrary) linear subspace
disjoint from $Y$. Assume  that $dim Q<r-m-1$ and that
$$2m>r+\sigma(Y)-1. \tag 2$$
Then the double locus of ${\pi_Q}_{|Y}$ (= locus of points $y\in
Y$ such that $\pi_Q\inv(\pi_Q(y))\cap Y\neq \{y\}$ as schemes) is
irreducible (resp. connected).\endproclaim Now the case where
$\dim Q>0$ (no other hypotheses needed but non-projectability) is
an easy and well known consequence, due to Franchetta and Mumford,
of Bertini's Theorem (see \cite{Mo}, p.115 or \cite{Pi} or below),
so the only new conclusion is when $Q$ is a point and as usual,
the case $Q$ arbitrary follows easily by connectedness principles
from the case $Q$ generic. For this case, it is convenient to
shift our viewpoint slightly, as follows.\par Let $Y^\2$ denote
the normalization of the blow-up of $Y\times Y$ along the diagonal
$\Delta_Y,$ with exceptional divisor $E_Y$. Let $I_Y$ denote the
tautological $\P^1$-bundle (or 'incidence variety') over $Y^\2$ (=
pullback of analogous object over $(\P^r)^{[2]}$). Being a $\P^1$
bundle over a normal variety,  $I_Y$ is also normal and we have a
diagram
$$\matrix\ I_Y&\overset f\to\to\P^r\\
\pi\downarrow&\\
\ Y^\2,&\endmatrix\tag 3$$ where non-projectability means $f$ is
surjective. Now it follows easily from the classical trisecant
Lemma (most secants are not multisecant) that for a generic linear
subspace $Q$, $f\inv(Q)$ is birational to the double locus of
${\pi_Q}_{|Y}$. If $\dim(Q)>0$, then $f\inv(Q)$ is automatically
irreducible by Bertini's theorem, which already proves Theorem 1
for this case. This result is due originally to Franchetta in the
case of surfaces (cf. [Fr], [En]); the foregoing argument is due
to Mumford and is given in [Mo].\par Let us say that $Y$ has the
{\it{irreducible secant cone}} (ISC) property if a generic fibre
of $f$ is irreducible. Then the remaining case $\dim(Q)=0$ of
Theorem 1 follows from (indeed, is equivalent to) the following
\proclaim{Theorem 2} Hypotheses as in Theorem 1, $Y$ has the ISC
property.\endproclaim In view of Zak's Tangency theorem, Theorem 2
implies the following result (which will shortly be improved
below): \proclaim{Corollary 3(temp)} If $Y\subset\P^r$ is smooth,
non-projectable with
$$ \dim(Y)>\frac{2}{3}(r-1),$$
 then $Y$ has the ISC property.\endproclaim
As another
 application, we obtain a proof of a version of Zak's linear normality
theorem (cf. [Z2], Thm. II.2.14): \proclaim{Corollary 4} Let
$X\subset\P^N$ be irreducible, nondegenerate, and set $b=0$ if $X$
is smooth and otherwise $b=\text{dim Sing}(X).$ Assume that
$$\dim(X)> \frac{1}{3}(2N+b-1).$$
Then $X$ is linearly normal, i.e. not the image of the bijective
projection of a nondegenerate subvariety of
$\P^{N+1}$.\endproclaim \demo{Proof that Corollary 3(temp) $=>$
Corollary 4} We use induction on $\dim(X)$. By generically
projecting, we may assume $X$ is nonprojectable. Assume for
contradiction that $X$ is the bijective projection of a
nondegenerate variety
$$\tilde{X}\subset\P^{N+1}$$
from a point $Q\in\P^{N+1}-\tilde{X}.$
 Let
$M\subset \P^N$ be a generic $\P^{N-b-1}$ and $Y=X\cap M$. Thus
$Y$ is smooth and spans $M$. Assume to begin with that $Y$ is
nonprojectable within $M$. By Corollary 3(temp), the secant
variety $\Sec(Y)$ coincides with $M$ 'with multiplicity 1', in the
sense that, for a generic linear $\P^{b+1}, L\subset \P^N$, the
scheme-theoretic inverse image $f\inv(L)$ is a reduced irreducible
subvariety of codimension $N-b-1$ lying over a single point of $L$
(viz. $L\cap M$); i.e. $f\inv(L)$ coincides with the fibre of
$I_Y\to \Sec(Y)=M$ over (the general point) $L\cap\Sec(Y)$. \par
On the other hand, note that $M$ is the projection of a unique
codimension-$(b+1)$ linear subspace of $\P^{N+1}$ containing $Q$,
say $A$, and $Y\simeq \tilde{X}\cap A$. Consequently,
 $Y$ can be viewed as a specialization of a
smooth subvariety
$$Y'\subset X,$$
which is the (isomorphic) projection of a generic
codimension-$(b+1)$ linear space section $\tilde{X}\cap
A'\subset\P^{N+1}$, with $Q\not\in A'$. Note that $Y'$ spans a
$\P^{N-b}$ which we denote by $M'$. By semi-continuity, similar
assertions as for $f$ must hold also for the analogously-defined
map $f_{Y'}$; thus $f_{Y'}\inv(L)$ is reduced, irreducible and of
codimension $N-b-1=\ $codim$(L)$. This implies firstly that
$\Sec(Y')$ is $(N-b-1)$-dimensional; then since $f_{Y'}\inv(L)$
has a component over each point of $L\cap\Sec(Y')$, the only way
$f_{Y'}\inv(L)$ can be irreducible is if $\Sec(Y')$ is a linear
$\P^{N-b-1}$, which contradicts the fact that $Y'$ spans $M'$ of
dimension $N-b$.\par Finally, suppose $Y$ is projectable within
$M$ and let $Y', M'$ be as above.
 Then the nondegenerate subvariety
 $$Y'\subset
M'=\P^{N-b}$$ of dimension $\dim(X)-b-1$ is projectable to
$\P^{N-b-2}$, which contradicts our induction hypothesis. \qed
\enddemo
\remark{Remarks} 1. For $X$ smooth, Zak's linear normality theorem
covers the larger range $\dim(X)>\frac{1}{3}(2N-2)$.\par 2. The
basic idea of the foregoing argument goes back to \cite{R}, and a
similar idea was recently used by Brandigi \cite{B} to prove
linear normality in the range $\dim(X)\geq\frac{3}{4}N.$
In fact, this argument proves the following fact of independent
interest: if $X\subset \P^N$ irreducible, nondegenerate and its
general hyperplane section is smooth, nonprojectable and has the ISC
property, then $X$ is linearly normal (in the above sense).
\par
3. Corollary 4 is sharp: to see this let $Z$ be a smooth Severi
variety (cf. [Z2]) in $\P^r$, that is $Z$ is projectable and
$\dim(Z)=\frac{2}{3}(r-2)$. Set $N = r+b$ and let $Z' \subset\P^{N+1}$
be the cone over $Z$ with vertex a $\P^b$, $X\subset \P^N$ its generic
(isomorphic) projection. Then $\dim(X)=\frac{1}{3}(2N+b-1)$ and
$\dim Sing(X)=b$.
\endremark
Given Corollary 4, we can sharpen slightly the statement of
Corollary 3(temp): \proclaim{Corollary 3} Let $Y\subset \P^r$ be
smooth nondegenerate with
$$\dim(Y)>\frac{2}{3}(r-1).$$
Then $Y$ is non-projectable and has the ISC property.\endproclaim
\demo{Proof} By Corollary 3(temp), it suffices to prove that $Y$
is non-projectable. If not, apply Corollary 4 to the generic
projection of $Y$ to $\P^{r-1}$ to deduce a
contradiction.\qed\enddemo
\remark {Remark} Again Corollary 3 is sharp: for this let
$Y\subset\P^r$ be the generic projection of a Severi variety
(cf. Remark 3 following Corollary 4). Then $Y$ is smooth,
non-projectable and does not have the ISC property
(e.g. because the cone on $Y$ is not linearly normal).
\endremark\par
It is amusing, perhaps, to translate
the irreducibility conclusion of Corollary 3 into cohomology
(taking for granted the nonprojectability conclusion).
Let $F$ denote a general fibre of $f$. Then $F$
is smooth, nonempty and $(2\dim(Y)+1-r)$-dimensional.
 Clearly irreducibility (i.e. connectedness)
of $F$ is equivalent, provided $q(Y)=h^1(\O_Y)=0,$ to the
vanishing
$$H^1(I_Y, \I_F)=0,\tag 4$$
where $\I_F$ denotes the ideal sheaf of $F$;
in the dimension range in question, $q(Y)=0$ automatically
by the Barth-Larsen Theorem. Pulling back the
Koszul resolution of the ideal sheaf of a point in $\P^r$ and
using standard vanishing results (e.g. [SS], Thm. 7.1) which imply
that
$$H^i(I_Y,f^*(\O(-j)))=0, \forall i<r, j>0,$$
 we see easily that (4)
is equivalent to the vanishing
$$H^r(I_Y, f^*(\O(-r)))=0.\tag 5$$
Let $\E$ denote the tautological subbundle on $Y^{[2]}$, so that
$I_Y=\P(\E).$ Then by standard computations the vanishing (5)
reduces to the vanishing on $Y^{[2]}$:
$$H^{r-1}(Y^{[2]}, \Sym^{r-2}(\E^v)\otimes \det(\E^v))=0.\tag 6$$
\proclaim{Corollary 5} With hypotheses as in Corollary 3,
 the vanishing (6) holds.
\endproclaim
Trying to find a direct proof of
 Corollary 5 might seem like a promising
route to a cohomological proof of Corollary 3, but we were unable
to find such a direct proof. This still looks like an intriguing,
though difficult problem.
\par
We now give the proof of Theorem 2, letting notations and
assumptions be as there. The basic idea is the following. Consider
a Stein factorization of $f$:
$$I_Y\to Z\overset g\to\to\P^r,$$
where $Z$ is normal and $g$ is generically finite and surjective.
Now it is a general fact that if $h:W\to T$ is a morphism
of irreducible varieties and $W$ is normal,
then so is a general fibre of $h$: this can
be seen, e.g. using Serre's criterion, or alternatively, use [G], Thm
12.2.4, which says, in the scheme-theoretic context, that the set
$N(h)$ of points $t\in T$ such that $h\inv(t)$ is normal is open;
when $W$ is normal, $N(h)$ contains the generic point (in the
scheme-theoretic sense) of $T$, hence also an open set of closed
points. In our case, since $I_Y$ is normal, it follows
 that so is a generic fibre of $f$,
therefore the irreducible and connected components of this fibre
coincide (cf. [E], Thm. 18.12). Consequently the degree of $g$
coincides with the number of irreducible (=connected) components
of the general fibre of $f$, so the Theorem's assertion is that
$g$ is birational.

 Then there is a Zariski open $U\subset\P^r$
such that $\P^r-U$ has codimension $>1$ and $g\inv(U)\to U$ is
finite, and we may assume $g\inv(U)$ is smooth as well. Since $U$,
like $\P^r$, is simply connected, it follows that if $\deg(g)>1$
then $g$, hence $f$ is {\it{ramified in codimension 1}}, i.e.
there is a prime divisor $F\subset I_Y$ such that $f(F)\subset
\P^r$ is a divisor and $f$ is ramified on $F$. We proceed to show
that the latter conclusion leads to a contradiction.\par Now it is
an easy consequence of the Fulton-Hansen Connectedness Theorem
(cf. \cite{FL}, Corollary 5.5) that in our case we have
$$f(\pi\inv(E_Y))=\P^r,$$
hence $F\neq\pi\inv(E_Y)$, and therefore a general point of $F$ is
of the form $(x,y,z)$ where $x,y\in Y$ are distinct and
$$z\in <x,y>$$
($<x,y>$ denotes the line spanned by $x,y$). Now a standard
computation known as Terracini's Lemma \cite{FR} says that
$$\im\ df_{(x,y,z)}=<T_xY,T_yY>,$$
and in particular this image is independent of $z\in <x,y>.$ It
follows that $F$ is the pullback of a divisor $D\subset Y^\2$,
where a general point $(x,y)\in D$ has the property that $x\neq y$
and
$$\rho:=\dim <T_xY,T_yY>\ <r.$$
%$$H_{x,y}:=<T_xY,T_yY>\subsetneq\P^r.\tag ??$$
%Now assume that $H_{x,y}$ is generically a hyperplane.
We may assume that the projection map $p_1:D\to Y$ is surjective,
and let $D_x\subset Y$ denote the image of its general fibre under
$p_2$, which is a divisor on $Y$. Setting $W=T_xY$, note that
 a general $y\in D_x$ is smooth on $Y$ and we have
$$\rho-1\leq\dim <T_yD_x, W>\leq \rho.\tag 4$$
Now consider the following diagram, with vertical arrows only
rational maps induced by projection from $W$:
$$\matrix S_{x,v}&\subset & D_x&\subset & \P^r&\\
:&&:&&:&\pi_W
\\
\downarrow&&\downarrow&&\downarrow\\
v&\in&V_x&\subset&\P^{c-1}_x.&\endmatrix
$$
Here $c=r-m, V_x$ is the (closure of the) image of $D_x$, $v\in
V_x$ is a general point and $S_{x,v}=\pi_W\inv(v)$, which we may
assume contains a general point of $D_x$. By (4), the dimension of
$V_x$ is either $\rho-m-1$ or $\rho-m-2$, and in these respective
cases we have $\dim S_{x,v}=2m-\rho$ (resp. $2m-\rho+1$). Though
not essential for our purposes, it is interesting to note that
when $x$ is viewed as variable, a general hyperplane in $\P^{c-1}$
corresponds in $\P^r$ to a general tangent hyperplane to $Y$, i.e.
a general element of the dual variety $Y^*.$\par Now suppose that
$V_x$ is of dimension $\rho-m-1$, so that $S_{x,v}$ is of
dimension $2m-\rho$. Note that by (4) this implies that for
general $y\in D_x,$
$$<T_yY,W>=<T_yD_x, W>, $$
which projects modulo $W$ to $T_vV_x, v=\pi_W(y)$. Now for a
linear subspace $U\subset \P^{c-1}_x,$ we denote by $\pi_W^*(U)$
the unique linear subspace of $\P^r$ which contains $W$ and
projects to $U$ (this is uniquely determined by $U$). Then we
conclude that for general $y\in S_{x,v},$ we have
$$T_yY\subset \pi_W^*(T_vV_x).$$
Thus the linear space $\pi_W^*(T_vV_x)$ of dimension $\rho\leq
r-1$ is tangent to $Y$ along a locus of dimension at least
$2m-\rho\geq 2m-r+1$, contradicting (2).\par Suppose now that
$V_x$ is of dimension $\rho-m-2$, so $S_{x,v}$ is of dimension
$2m-\rho+1$. Then for $y\in S_{x,v}$, the projection of $T_yY$ to
$\P^{c-1}_x$ is a $\P^{\rho-m-1}$ containing
$T_vV_x=\P^{\rho-m-2}$, and the set of all these linear subspaces
is a $\P^{r-\rho}$, so imposing such a subspace to stay fixed is
$r-\rho$ conditions. Thus, pulling back to $\P^r,$ we can find a
subvariety $T$ of codimension at most $r-\rho$ in $S_{x,v}$,
containing a general point of $S_{x,v}$ (hence of $D_x$), such
that $T_yY$ is contained in a fixed $\P^\rho$ for all $y\in T$.
 Since
$$\dim T\geq 2m-\rho+1-(r-\rho)=2m-r+1,$$
this again contradicts (2).\qed \remark{Example} By Corollary 3,
any smooth 3-fold in $\P^5$ has the ISC property. On the other
hand, if $Y$ is a smooth surface in $\P^4$, Theorem 2 says that
$Y$ has the ISC property {\it unless} $\sigma(Y)=1$, i.e. unless
$Y$ admits a hyperplane section with a multiple component. For
example, the projected Veronese surface  $Y\subset\P^4$ admits
multiple hyperplane sections, and indeed does not have the ISC
property; in fact, the double curve of its generic projection to
$\P^3$ consists of 3 conics on $Y$ mapping 2:1 to 3 lines on
$\bar{Y}\subset\P^3$. See \cite{GH}, pp. 628-635 for this and
other interesting examples. As we mentioned above, Franchetta
proved that for any smooth nondegenerate surface in $\P^5$ or
higher, other than the Veronese, the double curve of its generic projection to $\P^3$ is
irreducible.\endremark
\subheading{Acknowledgment} This work was begun when the second
author was visiting the Mathematics Department at Universit\`a di
Roma Tre. He would like to thank the department, especially
Edoardo Sernesi and Sandro Verra, the Director, for their
assistance and their hospitality.

%\bibliography
%\endbibliography
\enddocument